\documentclass{article}
\usepackage{arxiv}

\usepackage[utf8]{inputenc} 
\usepackage[T1]{fontenc}    
\usepackage{microtype}      
\usepackage{doi}

\usepackage{amsmath,amsfonts}
\usepackage{algorithmic}
\usepackage{algorithm}
\usepackage{array}
\usepackage{textcomp}
\usepackage{stfloats}
\usepackage{url}
\usepackage{verbatim}
\usepackage{graphicx}

\usepackage[style=numeric-comp,sorting=none,maxbibnames=10]{biblatex}
\usepackage{booktabs}
\usepackage{makecell}
\usepackage{hyperref}

\usepackage[capitalise]{cleveref}
\crefformat{equation}{#2(#1)#3}
\Crefformat{equation}{#2(#1)#3}
\crefrangeformat{equation}{#3(#1)#4--#5(#2)#6}
\Crefrangeformat{equation}{#3(#1)#4--#5(#2)#6}

\usepackage{csquotes}
\usepackage{siunitx}
\usepackage{mathtools}
\DeclareSIUnit{\euro}{\text{€}}
\DeclareSIUnit{\year}{\text{a}}
\DeclareSIUnit{\voltamperereactif}{\text{var}}

\title{The Benefits of an Integrated Approach for Stability-Constrained Power System Planning}
\date{}

\usepackage{authblk}

\newbox{\orcid}\sbox{\orcid}{\includegraphics[scale=0.06]{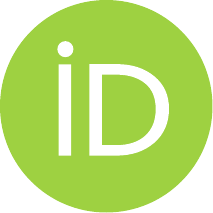}} 
\author[1,2]{%
	\href{https://orcid.org/0009-0000-6286-147X}{\usebox{\orcid}\hspace{1mm}Gereon~Recht\thanks{\texttt{gereon.recht@dlr.de}}}%
}
\author[2]{%
	\href{https://orcid.org/0009-0003-4504-5165}{\usebox{\orcid}\hspace{1mm}Benedikt~Jahn}%
}
\author[1,2]{%
	 \href{https://orcid.org/0000-0002-0080-7862}{\usebox{\orcid}\hspace{1mm}Oussama~Alaya}%
}
\author[1]{%
	\href{https://orcid.org/0000-0002-9720-0337}{\usebox{\orcid}\hspace{1mm}Karl-Kiên~Cao}%
}
\author[2]{%
	\href{https://orcid.org/0000-0002-0208-4100}{\usebox{\orcid}\hspace{1mm}Hendrik~Lens}%
}
\affil[1]{German Aerospace Center (DLR), Institute of Networked Energy Systems, Curiestr. 4, 70563 Stuttgart, Germany}
\affil[2]{Institute for Energy Process Engineering and Dynamics in Energy Systems (IED), University of Stuttgart, Germany}

\hypersetup{
pdftitle={The Benefits of an Integrated Approach for Stability-Constrained Power System Planning},
pdfauthor={Gereon~Recht, Benedikt~Jahn, Oussama~Alaya, Karl-Kiên~Cao, Hendrik~Lens},
}

\begin{document}

\maketitle

\begin{abstract}
Increasing penetration of inverter-based resources in today's power systems requires substitution of the contribution from synchronous generators to dynamic voltage stability and inertial response.
However, established approaches for power system planning are sequential, as stabilising measures are only considered at a later stage of the planning process.
We investigate the advantages of an integrated approach for power system planning, where stabilising measures are considered simultaneously with the expansion of generation, transmission, and storage systems via simplified stability constraints on inertia and voltage stability.
We find that system costs are reduced with the integrated approach and that the dual-use option of grid-forming battery energy storage systems is favoured over other stabilising measures like static synchronous compensators.
\end{abstract}

\section{Introduction}
With increasing penetration of inverter-based resources (IBRs), the contribution of synchronous generators (SGs) to inertial response and voltage stability must be substituted.
Yet power system planning, e.\,g.\ by transmission grid operators, today is still done by a sequential approach where stabilising measures are considered ex post \cite{systemstabilitaetsbericht}.
We investigate the advantages of an integrated approach, where stability is considered in power system planning simultaneously with expansion of generation, transmission, and storage systems via simplified constraints.
Least cost system designs with improved stability are valuable, as they can offer structural insight e.\,g.\ for grid connection regulations or siting of stabilising measures, and furthermore simplify the planning process.

Previous work includes \cite{10863371,11305559}, where rate of change of frequency (RoCoF) constraints for the centre of inertia (CoI) were considered, especially in the case of regional islanding or system splits.
We extend this work by including a constraint on dynamic voltage stability from \cite{10.1063/1.5002889} and comparing the integrated to the sequential approach.
In the experiment, we consider static synchronous compensators (STATCOMs).

\section{Methodology}
We extend the linear programming formulation of the \textsc{Capacity Expansion Problem} for aggregated network models described in \cite{11305559}.
The objective
\begin{align}
    & \min \sum_{s \in \mathcal{S}} \left( c_s u_s + \sum_{t \in \mathcal{T}} \left( o_s p_{s,t} + o_s^{\text{su}} \beta_{s,t}^{\text{su}} \right) \right) + \sum_{l \in \mathcal{L}} c_l u_l \label{eqn:obj}
\end{align}
represents the minimisation of annual system costs comprising capacity costs $c_s$, $c_l$ for expansion of power source units $u_s \in [u_s^{\text{min}}, u_s^{\text{max}}]$ and transmission line circuits $u_l \in [u_l^{\text{min}}, u_l^{\text{max}}]$ as well as operational costs $o_s$, $o_s^{\text{su}}$ of power injections $p_{s,t}$ and startup procedures $\beta_{s,t}^{\text{su}}$ over the set of snapshots $\mathcal{T} = \left\{ 1, \dots, T \right\}$.
The set of power sources $\mathcal{S}$ includes storage systems $\mathcal{S}^{\text{sto}}$ and power sources with grid-forming capabilities $\mathcal{S}^{\text{gfm}}$.
We write $\mathcal{S}_n$ to denote the set of power sources attached at bus $n$.
In the constraints
\begin{align}
    & 0 \leq \beta_{s,t} \leq u_s & \forall s \in \mathcal{S}, t \in \mathcal{T}, \\
    & \beta_{s,t} \leq \beta_{s,t-1} + \beta_{s,t}^{\text{su}}  & \forall s \in \mathcal{S}, t \in \mathcal{T} \setminus \{1\}, \\
    & \beta_{s,1} \leq \beta_{s,T} + \beta_{s,1}^{\text{su}}  & \forall s \in \mathcal{S}, \\
    & p_s^{\text{min}} \beta_{s,t} \leq p_{s,t} \leq a_{s,t} p_s^{\text{max}} \beta_{s,t} & \forall s \in \mathcal{S}, \forall t \in \mathcal{T},
\end{align}
we model the behaviour of power sources, where $\beta_{s,t}$ is the number of online units, $p_s^{\text{min}}$, $p_s^{\text{max}}$ are the power injection limits of one online unit, and $a_{s,t}$ is an exogenously given availability of supply.
Storage systems are power sources with additional variables and constraints
\begin{align}
    & 0 \leq p_{s,t}^c \leq a_{s,t} p_s^{\text{max}} u_s & \forall s \in \mathcal{S}^{\text{sto}}, \forall t \in \mathcal{T}, \\
    & p_{s,t} + p_{s,t}^c \leq a_{s,t} p_s^{\text{max}} \beta_{s, t}  & \forall s \in \mathcal{S}^{\text{sto}}, \label{eqn:complementarity_relaxation} \\
    & 0 \leq e_{s,t} \leq e_s^{\text{max}} u_s & \forall s \in \mathcal{S}^{\text{sto}}, \forall t \in \mathcal{T}, \\
    & e_{s,t} = e_{s,t-1} - \eta_s^{-1} p_{s,t} + \eta_s^c p_{s,t}^c + p_{s, t}^{\text{in}}   & \forall s \in \mathcal{S}^{\text{sto}}, \forall t \in \mathcal{T} \setminus \{1\}, \\
    & e_{s,1} = e_{s,T} & \forall s \in \mathcal{S}^{\text{sto}}.
\end{align}
Here, $p_{s,t}^c$ represents charging of the storage system, $e_{s,t}$ its state of charge, and $p_{s,t}^{\text{in}}$ a curtailable inflow.
Simultaneous charging and discharging is limited via \cref{eqn:complementarity_relaxation}.
Charging and discharging is associated with efficiencies $\eta_s^c$, $\eta_s$.
We differentiate between alternating current (AC) and direct current (DC) transmission lines $\mathcal{L}^{\text{ac}}, \mathcal{L}^{\text{dc}} \subseteq \mathcal{L}$.
For each AC line $l = (n,m)$, a positive power flow $p_{l,t}$ implies flow from bus $n$ to bus $m$, while a negative power flow implies flow in the opposite direction.
Moreover, losses are considered via the variable $p_{l,t}^{\text{loss}}$.
For DC lines, we define variables $p_{l,t}^{\text{al}}, p_{l,t}^{\text{ag}}$ for flow along and against the orientation of the line.
The power flow is then described by
\begin{align}
    & \sum_{l \in \mathcal{L}^{\text{ac}}} C_{lc} x_l p_{l,t} = 0 & \forall c \in \mathcal{C}, t \in \mathcal{T}, \\
    & p_{l,t}^{\text{loss}} \geq r_l p_{l0} (2p_{l,t} - p_{l0}) \begin{split}
        & \forall l \in \mathcal{L}^{\text{ac}}, t \in \mathcal{T}, \\
        & p_{l0} \in \mathcal{P}_{l0},
    \end{split} \label{eqn:losses} \\
    & \vert p_{l,t} \vert \leq a_l p_l^{\text{max}} u_l - p_{l,t}^{\text{loss}} & \forall l \in \mathcal{L}^{\text{ac}}, t \in \mathcal{T}, \\
    & 0 \leq p_{l,t}^{\text{al}}, p_{l,t}^{\text{ag}} \leq p_l^{\text{max}} u_l & \forall l \in \mathcal{L}^{\text{dc}}, t \in \mathcal{T}, \\
    \begin{split}
        &\sum_{s \in \mathcal{S}_n} p_{s,t} - \sum_{s \in \mathcal{S}_n^{\text{sto}}} p_{s,t}^c - d_{n,t} = \Biggl( \sum_{l \in \mathcal{L}^{\text{ac}}} K_{nl} p_{l,t} + \left\vert K_{nl} \right\vert \frac{p_{l,t}^{\text{loss}}}{2} \\
        &+ \sum_{\substack{l = (n,m)\\ \in \mathcal{L}^{{\text{dc}}}}} p_{l,t}^{\text{al}} - (1 - \eta_l) p_{l,t}^{\text{ag}} + \sum_{\substack{l = (m,n) \\ \in \mathcal{L}^{{\text{dc}}}}} p_{l,t}^{\text{ag}} - (1 - \eta_l) p_{l,t}^{\text{al}} \Biggr) \\
    \end{split} & \forall n \in \mathcal{N}, t \in \mathcal{T}, \label{eqn:kcl}
\end{align}
where $K_{nl}$, $C_{lc}$ are the incidence and cycle incidence matrices, $x_l$ and $r_l$ are the line reactance and resistance, $\mathcal{P}_{l0}$ is a set of points at which the tangents of the loss relaxation in \cref{eqn:losses} are centred, $p_l^{\text{max}}$ is the thermal limit of $l$, $d_{n,t}$ is the load, and $\eta_l$ is the loss factor of DC line $l$ \cite{11305559}.
To account for the change of impedances in transmission lines with a change of the number of parallel circuits, we solve the problem via sequential linear programming as in \cite{recht2026improvingoperationalfeasibilitylargescale}.
We moreover define the $\text{CO}_2$ limit
\begin{align}
    & \sum_{s \in \mathcal{S}, t \in \mathcal{T}} \frac{\xi_s}{\eta_s} p_{s,t} \leq \text{CO}_2^{\text{max}} \label{eqn:co2_limit}
\end{align}
with the $\text{CO}_2$ emissions per unit of thermal energy $\xi_s$ and the efficiency of the power source $\eta_s$.
We include three stability constraints
\begin{align} 
    &\Delta \dot{f}_a^{\text{max}} 2 \sum_{s \in \mathcal{S}_a} H_s p_s^{\text{max}} \beta_{s,t} \geq \vert \Delta P_a \vert & \forall a \in \mathcal{A}, t \in \mathcal{T}, \label{eqn:rocof_limit} \\
    & \Delta \dot{f}_a^{\text{max}} 2 \sum_{s \in \mathcal{S}_a} H_s p_s^{\text{max}} \beta_{s,t} \geq p_{a,t}^{\text{imb}} & \forall a \in \mathcal{A}^{\text{split}}, t \in \mathcal{T}, \label{eqn:rocof_limit_split} \\
    & \sum_{m \neq n} \vert B_{nm} \vert \leq \sum_{s \in \mathcal{S}_n^{\text{gfm}}} \frac{\beta_{s,t}}{x_s - x_s^\prime} & \forall n \in \mathcal{N}, t \in \mathcal{T}. \label{eqn:voltage_stability}
\end{align}
We limit the RoCoF of the CoI of each synchronous area $a \in \mathcal{A}$ via \cref{eqn:rocof_limit} with $\Delta \dot{f}^{\text{max}}_a = \dot{f}^{\text{max}}_a/f_a^b$, where $\dot{f}^{\text{max}}_a$ is the maximum RoCoF and $f_a^b$ the nominal frequency, the inertia constant $H_s$ and the predefined disturbance $\Delta P_a$.
In \cref{eqn:rocof_limit_split}, we limit the RoCoF for predefined system splits, defining the disturbance to be the power imbalance $p_{a,t}^{\text{imb}}$ due to power exchange before the separation \cite{11305559}.
We furthermore include \cref{eqn:voltage_stability}, a modified constraint from \cite{10.1063/1.5002889} on dynamic voltage stability, where $B$ is the nodal susceptance matrix and $x_s$, $x_s^\prime$ are the synchronous and transient reactances of the power source.
To the best of the author's knowledge, this constraint has not been included in a CEP formulation before.
In \cite{10.1063/1.5002889} a Kron-reduced network is assumed, i.\,e.\ an equivalent representation of the network reduced to the internal generator buses.
Since we consider aggregated representations of the power system, each electrical bus represents a larger region.
It is hence reasonable to assume that a Kron reduction is not necessary, since some grid-forming power source will always be operating at each bus.

\section{Numerical Experiment}
We conduct an experiment to compare the two approaches \textsc{Sequential} and \textsc{Stability} in two scenarios \emph{sg-dominated} and \emph{ibr-dominated} for the European power system.
In \textsc{Sequential}, an initial solution of \crefrange{eqn:obj}{eqn:co2_limit}, referred to as \textsc{Base} below, is reinforced with stabilising measures using the linear program defined by \cref{eqn:obj,eqn:rocof_limit,eqn:rocof_limit_split,eqn:voltage_stability}, where all relevant variables are fixed to the initial solution values except for stabilising measures, which only contribute to stability and do not affect normal operation.
\textsc{Stability} obtains a solution to the model \crefrange{eqn:obj}{eqn:voltage_stability}.
The scenarios differ with respect to the existing conventional generation in the problem instance and decarbonisation goals.
The scenario \emph{sg-dominated} contains all of today's existing conventional generation and aims to reduce $\text{CO}_2$ emissions to \qty{45}{\percent} of the levels of the year 1990.
In contrast, \emph{ibr-dominated} contains none of the existing generation with the exception of hydropower and aims at a full decarbonisation.
As a data basis, we use PyPSA-Eur (v2026.02.0) \cite{HORSCH2018207} at a resolution of 128 regions.
Stability-related assumptions adopted from \cite{11305559} include the disturbances $\Delta P_a$ for each synchronous area, which are determined either by the largest generator or HVDC connector, and a system split which separates the Iberian peninsula from the Continental Europe synchronous area.
We allow expansion of power sources, transmission, storage systems and stabilizing measures.
As stabilizing measures we consider STATCOMs \cite{nep_kosten} and battery energy storage systems with grid-forming inverters (BESS-GFM).
STATCOMs act as reactive power sources and can thus contribute to voltage stability.
E-STATCOMs are grid-forming STATCOMS that are additionally equipped with a supercapacitor in order to provide an inertial response with a typical inertia constant of \qty{6.25}{\second} in relation to the nominal power \cite{statcom_parameters,10045621}.
BESS-GFM are able to provide an inertial response at the expense of a reduced maximum power injection for normal operation, as the possible overcurrents of IBRs are lower than those of SGs.

\section{Results and Discussion}
\begin{figure}[t!]
    \centering
    \includegraphics{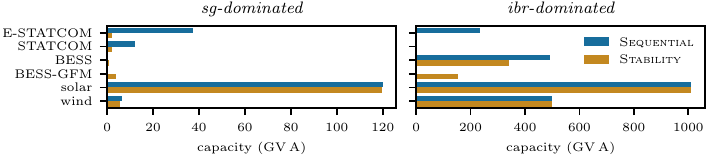}
    \caption{Generation and storage expansion results.}
    \label{fig:cap_mix}
\end{figure}

\cref{fig:cap_mix} shows that the strongest differences between the generation expansion results of \textsc{Stability} and \textsc{Sequential} in both scenarios are related to the choices of stabilising measures.
\textsc{Sequential} installs significantly more STATCOMs and E-STATCOMs, while \textsc{Stability} either keeps more SGs online in \emph{sg-dominated} or exclusively installs BESS-GFM in \emph{ibr-dominated}.
This leads to a reduction of system cost by \qty{0.27}{\percent} (\qty{1.024}{\giga\euro\per\year}) in \emph{sg-dominated} and \qty{2.28}{\percent} (\qty{6.507}{\giga\euro\per\year}) in \emph{ibr-dominated} when using \textsc{Stability} over \textsc{Sequential}.
While there are some local differences with respect to the expansion of solar and wind, there are only minor differences in the total capacities between the two models.
\begin{figure}[t!]
    \centering
    \includegraphics{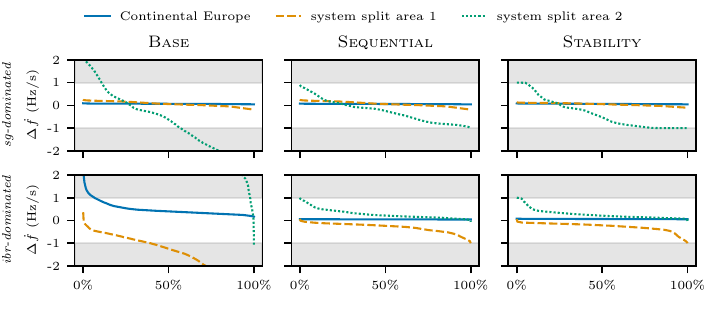}
    \caption{Duration curves for the CoI RoCoF of the Continental Europe synchronous area and both areas of the single system split considered.
    }
    \label{fig:rocof_duration_curves}
\end{figure}
\cref{fig:rocof_duration_curves} furthermore shows that the solution of \textsc{Base} contains many snapshots where RoCoF constraints are not fulfilled.
With \textsc{Base}, the voltage stability constraint \cref{eqn:voltage_stability} is violated for \qty{14.9}{\percent} of snapshots in \emph{sg-dominated} and for \qty{51.9}{\percent} in \emph{ibr-dominated}.
\textsc{Stability} shows more snapshots that are at the RoCoF limit than \textsc{Sequential}, suggesting that power exchange between the system split areas is reduced in favour of further installation of stabilising measures.
Solve times of \textsc{Stability} remain tractable, but increase in comparison to \textsc{Sequential} by \qty{65}{\percent} in \emph{sg-dominated} and \qty{68}{\percent} in \emph{ibr-dominated}.

One shortcoming of our experiment is the simplified representation of the sequential planning approach, where stability-related redispatch is ignored.
Moreover, we consider the CoI RoCoF, where larger local frequency deviations are ignored, and determine only a single system split.
While being a simplification, this is a common modeling shortcut in state-of-the-art planning processes \cite{systemstabilitaetsbericht}.
Lastly, the voltage stability constraint \cref{eqn:voltage_stability} is derived in \cite{10.1063/1.5002889} from a third-order model of synchronous generators in which voltage control of SGs is not considered.

\section{Conclusion}
We consider the advantages of an integrated approach for stability-constrained power system planning.
Our formulation includes constraints on voltage stability and the CoI RoCoF in the case of power injection deficits and system split disturbances.
In comparison to the sequential approach, the integrated approach leads to reduced system costs.
Its results furthermore show structural differences, as it favours the dual-use BESS-GFM over STATCOMs.
Moreover, we find that solutions which ignore stability generally show many violations of the stability constraints considered. 
Future research should investigate the inclusion of local RoCoF constraints and multiple system splits.

\section*{Acknowledgments}
The authors gratefully acknowledge funding from the German Federal Ministry for Economic Affairs and Energy under grant numbers FKZ 03EI1055A-B.

\printbibliography

\end{document}